# Special property of diagonal metrics of any dimensionality: an expression of the diagonal Riemann tensor components (RTC) in terms of the corresponding 2-d Gaussian curvature

Air1@nyu.edu

**Notation**

For a diagonal metric of any dimensionality (ie using orthogonal coordinates) we can write:

$$ds^2 = g_{\alpha\alpha}dx^\alpha dx^\alpha \quad \text{(summed on } \alpha\text{)}$$

**Notation**: We define: $g_D \equiv \sqrt{|g_{DD}|}$ ; **No summation over D**:
In the following, to indicate non-summation we will insert [**NS**].
We are not providing a geometric interpretation of this quantity (it is to be investigated) and do not claim it to be a tensor, and so the coordinate-letter D attached to g is not necessarily meant as an index, a 'subscript' or 'superscript'.

We also define $\eta_D$ and $\omega^D$ via:

$$\omega^D = g_D dx^D \quad ; \quad \textbf{NS: D} \quad \text{[ie: no summation over D]}$$

$$ds^2 \equiv \eta_\alpha (\omega^\alpha)^2 \equiv \eta_\alpha [g_\alpha]^2 [dx^\alpha]^2.$$

Where as can be seen the $\eta_\alpha$ contain the signature information, and we sum over $\alpha$. [Summation will be assumed in the below unless NS is specifically indicated.]   (1)   [1]

**The 2-d Gaussian curvature formula**

For a 2-d surface, using an orthogonal coordinate system $(x_A, x_B)$, the metric is:

$$ds^2 = g_{AA}dx^{A\,2} + g_{BB}dx^{B\,2} \quad \text{NS: A,B (ie no sum on A or B).}$$

Gauss's curvature formula for this surface is[2]:

$$K_{AB} = - [g_A g_B]^{-1} \{\eta_A [g_{A,B}/g_B]_{,B} + \eta_B [g_{B,A}/g_A]_{,A}\} \quad \text{NS: A,B.}$$

Note the manifest symmetry under exchange of A and B.

---

[1] **Note:** Since I define a one-index quantity $g_m = \text{sqrt}|g_{mm}|$, it is of interest to see that Landau Lifschitz (LL): p277: uses a similar notation for something else: $g_\alpha = g_{0\alpha}$, where $\alpha$ are the **spatial** indices only. [Also: tetrad notation in field theory, and tetrad formulation of GR in field theory.]

[2] See eg Barrett O'Neill: "Semi Riemannian Geometry" p81.

Terminology:
We'll say here that the Gaussian curvature arises from:
- the "interconnection of A & B", ie of the two dimensions represented by the coordinates corresponding to the two indices of K, or of the metric coefficients corresponding to those coordinates.
- a "mutual" interconnection, in the sense that 2nd derivatives arise via the metric coefficient corresponding to one coordinate derived w/r/t the other coordinate.

The above formula provides the curvature of 2-d surfaces in 3-d space, as well as the intrinsic curvature of a 2-d space described by the above metric.

**The formula for the Riemann Tensor Components (RTCs) for orthogonal metrics of any dimensionality:** We'll show below that it can be written as the Gauss term above plus an additional term which we call the 'intermediary term' I:

$$I^{AB} = \eta_D \eta_B \eta_A \; g_D^{-2} [g_{A,D}/g_A][g_{B,D}/g_B] \quad _D \neq _{A,B} \quad \text{(summed on D} \neq A,B)$$

$$= [g_A g_B]^{-1} \; \eta_D \eta_B \eta_A \; g_D^{-2} [g_{A,D}][g_{B,D}] \quad _D \neq _{A,B}$$

**Analysis**: This term vanishes unless *both* of the metric coefficients corresponding to the indices AB of the RTC have non-vanishing derivatives w/r/t *the same* 'other' (ie neither A nor B) coordinate(s.

**Terminology**: We will term these "intermediary coordinates" since they play an intermediary role in contributing non-zero curvature to the RTC. In other words, this term relates to the curvature arising from the interconnection of the two dimensions represented by A and B 'via' the dimensions represented by the "intermediary coordinates".

We can therefore write:  (Note: No Sum on A,B. However, D summed on D $\neq$ A,B)

$$R^{AB}{}_{AB} = K_{AB} + I^{AB} = -[g_A g_B]^{-1} \{ [g_{A,B}/g_B]_{,B} + \eta_A \eta_B [g_{B,A}/g_A]_{,A} + \eta_D \eta_B \eta_A \; g_D^{-2}[g_{A,D}][g_{B,D}] \; _D \neq _{A,B} \}$$

**Concluding analysis**: We can see that the RTC is composed of the Gaussian curvature arising from the "mutual" interconnection plus the curvature arising from the interconnection via the "intermediary" dimensions/coordinates.

- In contrast to the usual formulations of the RTC, the terms are manifestly symmetric under A $\leftrightarrow$ B. [(AR to AR: check the signature part of the last term)] Note that this is not true of the usual formulae for the RTC.
- This division of the RTC into "mutual" and "intermediary" terms grants this form of the RTC formula a certain intuitive basis.
- Use of this formula facilitates intuiting shortcuts in calculations of the RTC. (Some RTCs can be computed mostly 'by inspection' using this formula; the terms arising in the RTCs can be easily traced to the metric. Some examples of computation using the above formula will be provided. Some RTC properties of special metrics are easily discernible eg for metrics with inverse relationship eg: $g_{tt} = -g_{rr}^{-1}$ [eg $R^t{}_t$ is a Laplacian

# Derivation of the above RTC Formula

**Introduction**

We wish to find the Riemann tensor components (RTCs) $R^A{}_{BAB}$ from the formula for the curvature two-form:

$$\mathcal{R}^A{}_B = d\omega^A{}_B - \omega^A{}_D \wedge \omega^D{}_B \qquad [3] \qquad (2a)$$

$$= - R^A{}_{B|EF|}\, \omega^E \wedge \omega^F \qquad (E<F) \qquad [4]$$

From the above formula one can see that the (A,B,E,F) RTCs are obtained as the coefficients of $\omega^E \wedge \omega^F$.

Note that the RTC's are defined in this way, with only one index raised.

With both indices raised:

$$\mathcal{R}^{AB} = - R^{AB}{}_{|EF|}\, \omega^E \wedge \omega^F \qquad [5] \qquad (2b)$$

$$= - R^{AB}{}_{AB}\, \omega^A \wedge \omega^B \;_{(NS:A,B)} - R^{AB}{}_{AD}\, \omega^A \wedge \omega^D \;_{(D>A)} - R^{AB}{}_{DB}\, \omega^D \wedge \omega^B \;_{(D<B)} - R^{AB}{}_{CE}\, \omega^C \wedge \omega^E \;_{(C,E \neq A,B)}$$

For diagonal metrics, we will be interested in calculating the RTCs with repeating indices, ie the $R^A{}_{BAB}$ (NS: A,B) [6] and so we seek the coefficient of the term $\omega^A \wedge \omega^B$. We will however also inter alia obtain formulae for all the RTC components.

We begin below with a calculation of $\omega^A{}_B$ from which we will compute both terms in eq (2a); we'll refer to $d\omega^A{}_B$ as "term #1", and to $\omega^A{}_D \wedge \omega^D{}_B$ as "term #2".

**A) Computation of $\omega^A{}_D$:**

For the specific coordinate $x_A$: [A:NS]

$$d\omega^A = d[g_A dx^A] \qquad [A:NS]$$

$$[7] \qquad = d[g_A]dx^A + g_A \cancel{ddx^A} = g_{A,D}\, dx^D \wedge dx^A \qquad {}_D \neq {}_A$$

---

[3] see eg MTW p351
[4] see eg MTW p352
[5] see eg MTW p358
[6] **AR to AR:** show that only these exist, or are needed to set up or solve the Einstein equation.
[7] Note that by definition, operating twice with d, ie dd, gives zero.

$$= g_{A,D} \, \omega^D/g_D \wedge \omega^A/g_A = [g_{A,D}/(g_D g_A)] \, \omega^D \wedge \omega^A \quad {}_D \neq {}_A. \qquad (3)$$

However we also know that:

$$d\omega^A = -\omega^A{}_D \wedge \omega^D \quad {}_D \neq {}_A = +\omega^D \wedge \omega^A{}_D \quad {}_D \neq {}_A \qquad {}^8 \qquad (4)$$

By equating (3) and (4) we obtain:

$$\omega^D \wedge \omega^A{}_D \quad {}_D \neq {}_A = \omega^D \wedge \omega^A \, [g_{A,D}/(g_D g_A)] \quad {}_D \neq {}_A \qquad (5)$$

The first impression is that

$$\omega^A{}_D = \omega^A \, [g_{A,D}/(g_D g_A)] \quad {}_D \neq {}_A$$

however since by definition

$$\omega^D \wedge \omega^D = 0$$

adding a term proportional to $\omega^D$ has no effect, ie:

$$[\omega^A + f \, \omega^D] \wedge \omega^D = \omega^A \wedge \omega^D$$

where f is a function of the coordinates; we'll signify it's unknown nature by writing it below as [?].

Thus the more general solution of equation (5) above is:

$$\omega^A{}_D = [g_{A,D}/(g_D g_A)] \, \omega^A + [?] \, \omega^D \qquad (5b)$$

where the unknown function [?] will be determined below.

Note that whereas A is a specific coordinate, in contrast D is summed over all coordinates other than A, and so for each coordinate value of D there will be a companion term: $[?] \, \omega^D$.

We now interchange A ←→ D to obtain:

$$\omega^D{}_A = [g_{D,A}/(g_A g_D)] \, \omega^D + [?] \, \omega^A \qquad (5c)$$

We also know that:

$$\omega^A{}_D = -\omega^D{}_A \, \eta_A \, \eta_D \quad {}_D \neq {}_A$$

We will assume from this point on that we sum over all values of D except for A.

Therefore:

$$\omega^A{}_D = [g_{A,D}/(g_D g_A)] \, \omega^A + [?] \, \omega^D$$

$$= -\omega^D{}_A \, \eta_A \, \eta_D$$

---

[8] If D = A: $\omega^A{}_A = -\omega_A{}^A = -\eta_A \, \eta_A \, \omega^A{}_A = -\eta_A{}^2 \, \omega^A{}_A = -\omega^A{}_A$, so we do not allow D = A.

$$= - \eta_A \eta_D \{ [g_{D,A}/(g_A g_D)] \omega^D + [?] \omega^A \}$$

So we can set:

$$[?] \omega^D = - \eta_A \eta_D [g_{D,A}/(g_A g_D)] \omega^D \qquad (5d)$$

$$- \eta_A \eta_D [?] \omega^A = [g_{A,D}/(g_D g_A)] \omega^A \qquad {}^9$$

Therefore inserting $[?] \omega^D$ from (5d) into (5b) we can write the connection two-form as:

$$\omega^A{}_D = [g_{A,D}/(g_D g_A)] \omega^A - \eta_A \eta_D [g_{D,A}/g_A g_D)] \omega^D \qquad \text{Footnote 10}$$

$$= [g_A g_D]^{-1} \{ g_{A,D} \omega^A - \eta_A \eta_D g_{D,A} \omega^D \} . \qquad (6A1)$$

Written in terms of dx:

$$\omega^A{}_D = [g_A g_D]^{-1} \{ g_{A,D} g_A dx^A - \eta_A \eta_D g_{D,A} g_D dx^D \} \qquad (6A2)$$

$$= [g_{A,D}/g_D] dx^A - \eta_A \eta_D [g_{D,A}/g_A] dx^D . \qquad (6c)$$

Note the symmetry of the two terms in A and D (other than for the signature symbols).

## Significance of the above equation:

- Note from (6A1) that if both $g_{A,D}$ and $g_{D,A}$ are zero then $\omega^A{}_D = 0$. *This makes geometric sense: if both of these derivatives vanish, there is clearly no direct 'connection' between the coordinates A and D. Often we will find that the two metric coefficients $g_A$ and $g_D$ are functions of only one of the coordinates (A,D), not both, and so by inspection we will know that one of the two terms in (6) vanishes.*

- **The above formulation allows use of a given form of $\omega^A{}_D$ for all further computations!**

- **The development up to this point displays in a transparent manner the origin of each term in $\omega^A{}_D$.**

………………………….

**Defining the "sqrt connection":**

The connection component in terms of sqrt g's is:

---

[9] Interchanging D and A transforms one of the above two equations into the other into the other except that the signature factor is on the wrong side of the equaiotn. However since $-\eta_A \eta_D$ is simply a sign, it can be transferred from the rhs to the lhs without change. As a result one can easily see that these two equations are indeed the same.

[10] Remember that D will be summed over, and each $\omega^D$ will have the appropriate sign due to the accompanying $\eta_D$.

$$\Gamma^D{}_{AA} = (-1/2)\, g_{AA,D}/g_{DD}$$

$$= (-1/2)\, [\eta_A g_A{}^2]_{,D}/[\eta_D g_D{}^2]$$

$$= (2/-2)\, \eta_A g_A g_{A,D}/[\eta_D g_D{}^2]$$

$$= -\eta_D \eta_A\, [g_A/g_D][g_{A,D}/g_D]$$

Note the factor $g_{A,D}/g_D$ which displays a similar relationship between the metric components as in $g_{AA,D}/g_{DD}$. If we define the symbol:

$$\sqrt{\Gamma^D{}_{AA}} \equiv g_{A,D}/g_D$$

Which we will call "the sqrt connection" (which is NOT the sqrt of the connection), then we can write the connection component as:

$$\Gamma^D{}_{AA} = -\eta_D \eta_A\, [g_A/g_D]\, \sqrt{\Gamma^D{}_{AA}} \qquad . \qquad [11]$$

…………………………………………..

## B) Computation of $\omega^A{}_D \wedge \omega^D{}_B$

From equation (6A1)

$$\omega^A{}_D = [g_A g_D]^{-1}\{g_{A,D}\, \omega^A - \eta_A \eta_D\, g_{D,A}\, \omega^D\} \qquad (NS: A)$$

we can write:

$$\omega^A{}_D \wedge \omega^D{}_B \quad [D \neq A,B]$$

$$= [g_A g_D]^{-1}\{g_{A,D}\, \omega^A - \eta_A \eta_D\, g_{D,A}\, \omega^D\} \wedge [g_B g_D]^{-1}\{g_{D,B}\, \omega^D - \eta_D \eta_B\, g_{B,D}\, \omega^B\} \quad (NS: A,B) \qquad \textbf{(7a)}$$

[As a check: $\omega^A{}_D$ and $\omega^D{}_B$ have the same form, but with $(A,D) \rightarrow (D,B)$ ie $[A \rightarrow D,\ D \rightarrow B]$.

{Note: We will be interested a little later on in the cross-product containing $\omega^A \wedge \omega^B$, ie

$$[g_{A,D}\, \omega^A]\, [-\eta_D \eta_B\, g_{B,D}\, \omega^B] \quad \text{or} \quad -\eta_D \eta_B\, g_{A,D}\, g_{B,D}\, [\omega^A \wedge \omega^B] ,$$

---

[11] Note also that (see MTW p354): $\omega_{AD} \equiv g_{AN} \omega^N{}_D = \Gamma_{ADN} \omega^N$.

and note that this has signature symbols $\eta_D$ and $\eta_B$ but not $\eta_A$, so this term is NOT manifestly symmetric, which arises from the manifest non-symmetry between A and B in $\omega^A{}_D \wedge \omega^D{}_B$ .}

We combine the two products $[g_A g_D]^{-1}$ and $[g_B g_D]^{-1}$:

$$[g_A g_D]^{-1} [g_B g_D]^{-1} = [g_A g_D g_B g_D]^{-1} = [g_A g_B g_D^2]^{-1}$$

so that we can write:

$$\omega^A{}_D \wedge \omega^D{}_B \;_{[D \neq A,B]} = [g_A g_B g_D^2]^{-1} \{g_{A,D}\, \omega^A - \eta_A \eta_D\, g_{D,A}\, \omega^D\} \wedge \{g_{D,B}\, \omega^D - \eta_D \eta_B\, g_{B,D}\, \omega^B\}$$

Term:  #1   #2   #3   #4

Multiplying the terms in eq(7a) will give us these products:
$$(1) \times (3) + (1) \times (4) + (2) \times (3) + (2) \times (4)$$

Note that:
- $(2) \times (3) = 0$ because both have $\omega^D$ ;
- in $(2) \times (4)$ one has $- \times - = +$, and $\eta_D \eta_D = 1$, leaving only $\eta_A \eta_B$;

Therefore:
$$\omega^A{}_D \wedge \omega^D{}_B \quad D \neq A,B \quad (\text{NS: A,B}) = (1) \times (3) + (1) \times (4) + (2) \times (4)$$

$$= [g_A g_B g_D^2]^{-1}\{g_{A,D}\, g_{D,B}\, \omega^A \wedge \omega^D - \eta_D \eta_B\, g_{A,D}\, g_{B,D}\, \omega^A \wedge \omega^B + \eta_A \eta_B\, g_{D,A}\, g_{B,D}\, \omega^D \wedge \omega^B\}_{D \neq A,B} \quad \textbf{(7b)}$$

Note that the above is term #2 in eq (2a) .   {(eq 2a): $\mathcal{R}^A{}_B = d\omega^A{}_B - \omega^A{}_D \wedge \omega^D{}_B$ }

### C) Computation of $d\,\omega^A{}_B$

$$\omega^A{}_B = [g_{A,B}/g_B][dx^A] - \eta_A \eta_B [g_{B,A}/g_A][dx^B] \quad . \quad \text{(eq 6c with: D} \rightarrow \text{B)}$$

From equation 6c above we find:

$$d\,\omega^A{}_B = [g_{A,B}/g_B]_{,c}\, dx^C \wedge dx^A \quad C \neq A \quad - \eta_A \eta_B [g_{B,A}/g_A]_{,c}\, dx^C \wedge dx^B \quad C \neq B \qquad (8a) \quad (\text{NS: A,B})$$

Changing the order of the first wedge product, and therefore changing the sign:

$$d\,\omega^A{}_B = - [g_{A,B}/g_B]_{,c}\, dx^A \wedge dx^C \quad C \neq A \quad - \eta_A \eta_B [g_{B,A}/g_A]_{,c}\, dx^C \wedge dx^B \quad C \neq B \qquad (8a) \quad (\text{NS: A,B})$$

Note that in order to take the derivatives we employed eq 6C with $dx^A$ rather than 6A with $\omega^A$; however now that we have the derivative we can put this back into terms of $\omega^A$ and $\omega^B$:

$$d\,\omega^A{}_B = -[g_{AC}]^{-1}[g_{A,B}/g_B]_{,c}\,\omega^A \wedge \omega^C \quad C \neq A \quad - \eta_A \eta_B [g_C g_B]^{-1}[g_{B,A}/g_A]_{,c}\,\omega^C \wedge \omega^B \quad C \neq B \qquad (8b)$$

<div style="text-align:center">Term #1                Term #2</div>

Note:

- In the first term's sum over C, C cannot be A, but since B is NOT A, C CAN be B.
- when in the first term C = B, and in the second term C = A both terms contain $\omega^A \wedge \omega^B$.
- In the above case, $[g_B g_A]^{-1}$ is a common factor.

For term #1: We separate out the case C = B and C $\neq$ B.
For term #2: we separate below the cases C = A, and C $\neq$ A.           (NS: A,B)

$$d\,\omega^A{}_B = \text{a) Term #1 with } C = B \;\rightarrow\; -[g_A g_B]^{-1}\,[g_{A,B}/g_B]_{,B}\,\omega^A \wedge \omega^B$$

$$\text{b) + Term #1 with } C \neq B \;\rightarrow\; -[g_A g_C]^{-1}[g_{A,B}/g_B]_{,c}\,\omega^A \wedge \omega^C \quad C \neq A,B$$

$$\text{c) + Term #2 with } C = A \;\rightarrow\; -\eta_A \eta_B\,[g_A g_B]^{-1}\,[g_{B,A}/g_A]_{,A}\,\omega^A \wedge \omega^B$$

$$\text{d) + Term #2 with } C \neq A \;\rightarrow\; -\eta_A \eta_B[g_C g_B]^{-1}[g_{B,A}/g_A]_{,c}\,\omega^C \wedge \omega^B \quad C \neq A,B$$

Collecting terms with the same wedge products:

a) and c) have the same wedge: also, $[g_A g_B]^{-1}$ is a common factor: so:

$$\text{a) + c) } = -[g_A g_B]^{-1}\,\{[g_{A,B}/g_B]_{,B} + \eta_A \eta_B\,[g_{B,A}/g_A]_{,A}\}\,\omega^A \wedge \omega^B$$

In d), switching $\omega^C \wedge \omega^B$ to $\omega^B \wedge \omega^C$ so that we need to introduce a '-', and then adding to b):

$$\text{b) + d) } = -g_C^{-1}\,\{g_A^{-1}[g_{A,B}/g_B]_{,c}\,\omega^A - \eta_A \eta_B\,[g_B]^{-1}[g_{B,A}/g_A]_{,c}\,\omega^B\} \wedge \omega^C \quad C \neq A,B$$

and so we can write: note: we change the summation index from C to D:

$$d\,\omega^A{}_B = -[g_A g_B]^{-1}\,\{[g_{A,B}/g_B]_{,B} + \eta_A \eta_B\,[g_{B,A}/g_A]_{,A}\}\,\omega^A \wedge \omega^B +$$

$$-\{[g_A g_D]^{-1}[g_{A,B}/g_B]_{,D}\,\omega^A - \eta_A \eta_B[g_D g_B]^{-1}[g_{B,A}/g_A]_{,D}\,\omega^B\} \wedge \omega^D \quad D \neq A,B \quad (NS: A,B) \qquad (8c)$$

Note that the above is term #1 in eq (2a).    {eq (2a): $\mathcal{R}^A{}_B = d\omega^A{}_B - \omega^A{}_D \wedge \omega^D{}_B$ }

### D) Finding the Riemann tensor components (RTCs)

Inserting the results of equations (7b) and (8c) into eq. (2): (we use the summation index D):(NS: A,B) To more easily spot the relevant terms we'll use bold for $\boldsymbol{\omega^A \wedge \omega^B}$.

$$\mathcal{R}^A{}_B = d\omega^A{}_B - \omega^A{}_D \wedge \omega^D{}_B = - R^A{}_{BAB}\,\omega^A \wedge \omega^B - R^A{}_{BAD}\,\omega^A \wedge \omega^D - R^A{}_{BDB}\,\omega^D \wedge \omega^B$$

$$= -[g_A\,g_B]^{-1}\{[g_{A,B}/g_B]_{,B} + \eta_A\,\eta_B[g_{B,A}/g_A]_{,A}\}\,\omega^A \wedge \omega^B$$

$$\{-[g_A\,g_D]^{-1}[g_{A,B}/g_B]_{,D}\,\omega^A + \eta_A\,\eta_B[g_D\,g_B]^{-1}[g_{B,A}/g_A]_{,D}\,\omega^B\} \wedge \omega^D \quad D \neq A,B$$

$$- [g_A g_B g_D^2]^{-1}\{g_{A,D}\,g_{D,B}\,\omega^A \wedge \omega^D - \eta_D\eta_B\,g_{A,D}g_{B,D}\,\omega^A \wedge \omega^B + \eta_A\eta_B\,g_{D,A}\,g_{B,D}\,\omega^D \wedge \omega^B\}_D \neq A,B$$

We now collect terms with the same wedge products:

Note the three terms in $\omega^A \wedge \omega^B$: (the first two, and the second-to-last): adding them gives:

$$-[g_A\,g_B]^{-1}\{[g_{A,B}/g_B]_{,B} + \eta_A\,\eta_B[g_{B,A}/g_A]_{,A}\}\omega^A \wedge \omega^B - \{[g_A g_B g_D^2]^{-1}(-)\eta_D\eta_B\,g_{A,D}g_{B,D}\,\omega^A \wedge \omega^B\}_D \neq A,B$$

Distributing the common $\omega^A \wedge \omega^B$ and $-[g_A\,g_B]^{-1}$ gives us:

$$- [g_A\,g_B]^{-1}\{[g_{A,B}/g_B]_{,B} + \eta_A\,\eta_B[g_{B,A}/g_A]_{,A} - g_D^{-2}\,\eta_D\eta_B\,g_{A,D}g_{B,D}\}\omega^A \wedge \omega^B \quad D \neq A,B$$

Note that there are only two terms in $\omega^B \wedge \omega^D$:

$$1)\; - \eta_A\,\eta_B[g_D\,g_B]^{-1}[g_{B,A}/g_A]_{,D}\,\omega^B \wedge \omega^D \quad D \neq A,B$$

$$2)\; - [g_A\,g_B\,g_D^2]^{-1}\eta_A\eta_B\,g_{D,A}\,g_{B,D}\,\omega^D \wedge \omega^B \quad D \neq A,B$$

$$= \eta_A\eta_B[g_D\,g_B\,g_A\,g_D]^{-1}\,g_{D,A}\,g_{B,D}\,\omega^B \wedge \omega^D \quad D \neq A,B$$

Terms 1) and 2) have $\eta_A\eta_B[g_D g_B]^{-1}\omega^B \wedge \omega^D$ in common, and so adding 2) + 1) gives:

$$\eta_A\eta_B[g_D g_B]^{-1}\{[g_A g_D]^{-1}\,g_{D,A}\,g_{B,D} - [g_{B,A}/g_A]_{,D}\}\;\omega^B \wedge \omega^D$$

so we find:

$$\mathcal{R}^A{}_B = - [g_A\,g_B]^{-1}\{[g_{A,B}/g_B]_{,B} + \eta_A\,\eta_B[g_{B,A}/g_A]_{,A} - g_D^{-2}\,\eta_D\eta_B\,g_{A,D}\,g_{B,D}\}\omega^A \wedge \omega^B \quad D \neq A,B$$

$$+ \eta_A\eta_B[g_D g_B]^{-1}\{[g_A g_D]^{-1}\,g_{D,A}\,g_{B,D} + [g_{B,A}/g_A]_{,D}\}\;\omega^B \wedge \omega^D$$

$$- g_A^{-1}\{g_D^{-1}[g_{A,B}/g_B]_{,D} + [g_B g_D^2]^{-1}\,g_{A,D}\,g_{D,B}\}\;\omega^A \wedge \omega^D \quad D \neq A,B$$

$$= - R^A{}_{BAB}\,\omega^A \wedge \omega^B - R^A{}_{BAD}\,\omega^A \wedge \omega^D - R^A{}_{BDB}\,\omega^D \wedge \omega^B$$

Therefore:

$$R^A{}_{BAB} = + [g_A g_B]^{-1} \{ [g_{A,B}/g_B]_{,B} + \eta_A \eta_B [g_{B,A}/g_A]_{,A} - g_D{}^{-2} \eta_D \eta_B g_{A,D} g_{B,D} \}$$

$$= + [g_A g_B]^{-1} \{[g_{A,B}/g_B]_{,B} + \eta_A \eta_B [g_{B,A}/g_A]_{,A}\} - \eta_D \eta_B g_D{}^{-2} [g_{A,D}/g_A][g_{B,D}/g_B] \quad _D \neq _{A,B}$$

(check signs, and raised/lowered indices):

Note that as pointed out earlier the third term seems non-symmetric since only $\eta_D \eta_B$ appear without $\eta_A$. However we can see that this must be the case since the $\eta$'s appear in pairs, their presence makes a difference only if one of A or B is t, so that the product is – (whatever the choice of signature). If one or the other is t, then $\eta_A \eta_B$ is necessarily –1 and the sum of the first two terms become the negative of its sum under interchange A → B. If the third term had $\eta_A \eta_B$ in it then this third term would not change sign under A → B, so therefore it must have only one of them, and has one sign if A is t and another sign if B is t. ie when there are two terms involved as in the first part of $R^{AB}{}_{AB}$, the switch in sign can be overall by both terms conspiring to switch together, but for one term {the $[g_{A,D}/g_A][g_{B,D}/g_B]$ term} the sign switch can happen only if there's an $\eta$ for only one of the two coordinates.

…………

Note:
- The first term is exactly Gauss's curvature formula for a 2-d surface.
- For an orthogonal metric, the RTC's indices involve only two coordinates, and the RTC itself involves only the metric coefficients corresponding to these two indices, with the other coordinates appearing only in derivatives wrt them.

……………

**Note:**
Similarly, the non-diagonal components of the RTCs are (check signs, and raised/lowered indices):

$$R^{AB}{}_{BD} = - \eta_A \eta_B [g_D g_B]^{-1} \{ [g_A g_D]^{-1} g_{D,A} g_{B,D} + [g_{B,A}/g_A]_{,D} \}$$

$$R^{AB}{}_{AD} = - g_A{}^{-1} \{ g_D{}^{-1} [g_{A,B}/g_B]_{,D} + [g_B g_D{}^2]^{-1} g_{A,D} g_{D,B} \} \; \omega^A \wedge \omega^D \quad _D \neq _{A,B}$$

………………

…….

**Question**: Is there a way of characterizing invariantly those Riemann spaces (ie with a connection & metric) which are amenable to description in terms of a diagonal metric?
Perhaps there is a hint in this fact: all 2-d metric spaces can be written using a diagonal metric, and those higher-d spaces which can also be written using diagonal metrics are the 2-d curvature + an additional term involving the additional dimensions.
Associated to this question: what about a space with a connection but not a metric, and which is amenable to description in terms of orthogonal coordinates.
What is the geometric nature of the "sqrt connections" and of the single-index g's, like vierbeins etc
…
Sample calculations will be provided, for 4 and 5 dimensions.

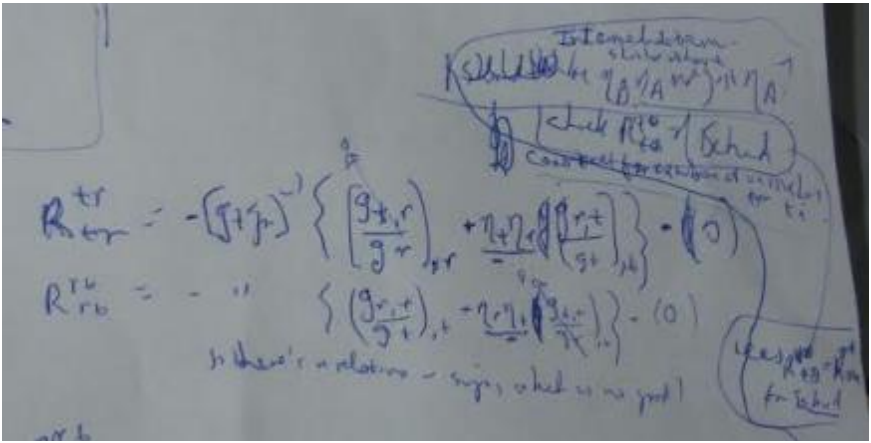

.....

..................

# Appendix
# Post-Hoc Heuristic 'Derivation" of the repeated-index RTC

**A) Finding $\omega^A{}_D$:**

From the way that the switching of upper to lower indices of $\omega^A{}_D$ changes the sign as follows:

$$\omega^A{}_D = -\omega^D{}_A \, \eta_A \, \eta_D {}_{.D} \neq {}_A$$

we can justify the following structure for $\omega^A{}_D$:

$$\omega^A{}_D = [g_{A,D}/(g_D g_A)] \, \omega^A - \eta_A \, \eta_D \, [g_{D,A}/g_A g_D)] \, \omega^D$$

$$= [g_A g_D]^{-1} \{g_{A,D} \, \omega^A - \eta_A \, \eta_D \, g_{D,A} \, \omega^D\} \quad . \tag{6A1}$$

Written in terms of dx: ie: $\omega^A = g_A dx^A$ ; $\omega^D = g_D dx^D$ :

$$\omega^A{}_D = [g_A g_D]^{-1} \{g_{A,D} \, g_A dx^A - \eta_A \, \eta_D \, g_{D,A} \, g_D \, dx^D\} \tag{6A2}$$

$$\omega^A{}_D = [g_{A,D}/g_D] \, [dx^A] - \eta_A \, \eta_D \, [g_{D,A}/g_A][dx^D] \quad . \tag{6c}$$

**B) Computation of the 'wedge product' $\omega^A{}_D \wedge \omega^D{}_B$**

To form $\omega^D{}_B$ we simply transform $\omega^A{}_D$ above by switching (A,D → D,B) ie (A→D), (D→B), so that the first term is in $\omega^D$ and the second has $\omega^B$. Creating the 'wedge product' of the two will lead to only one term with $\omega^A \wedge \omega^B$, basically from the product of the first term of $\omega^A{}_D$ (which has $\omega^A$) with the second term of $\omega^D{}_B$ (which has $\omega^B$); ie

$$[g_{A,D}\ \omega^A] \wedge [-\eta_D\ \eta_D\ g_{B,D}\ \omega^B] = -\eta_D\ \eta_B\ g_{A,D}\ g_{B,D}\ \omega^A \wedge \omega^B$$

with a factor of $[gg]^{-1}$ for both. Thus the relevant part $\omega^A{}_D \wedge \omega^D{}_B$ of is:

$$-\eta_D\ \eta_B\ g_{A,D}\ g_{B,D}\ [g_A g_D]^{-1}[g_D g_B]^{-1}\ \omega^A \wedge \omega^B \quad D \ne A,B$$

$$= -\eta_D \eta_B\ g_D^{-2}\ [g_{A,D}/g_A][g_{B,D}/g_B]\ _D \ne _{A,B} \qquad (7e)$$

C) $d\ \omega^A{}_B$ :

$$d\ \omega^A{}_B = -[g_{A,B}/g_B]_{,c}\ dx^A \wedge dx^C\ _{C \ne A} - \eta_A\ \eta_B\ [g_{B,A}/g_A]_{,c}\ dx^B \wedge dx^C\ _{C \ne B} \qquad (8a)\ (NS: A,B)$$

$$= -[g_{A,B}/g_B]_{,c}\ [\omega^A/g_A] \wedge [\omega^C/g_C]\ _{C \ne A} - \eta_A\ \eta_B\ [g_{B,A}/g_A]_{,c}\ [\omega^B/g_B] \wedge [\omega^C/g_C]\ _{C \ne B}$$

The terms of interest are for C = B in the first term, and for C = A in the second, which makes $\omega^A \wedge \omega^B$ a common factor, and turns the $g_i$'s under the $\omega^i$'s, ie both $g_A g_C$ and $g_B g_C$ into $g_A g_B$ and thereby makes $[g_A g_B]^{-1}$, a common factor as well. Distributing these common factors we have:

$$d\ \omega^A{}_B = -[g_A\ g_B]^{-1}\ \{[g_{A,B}/g_B]_{,B} - \eta_A\ \eta_B\ [g_{B,A}/g_A]_{,A}\}\omega^A \wedge \omega^B \qquad (8c)\ (NS: A,B)$$

Note that the above is the (ABAB) part of term #1 in eq (2a).    [eq (2a): $\mathcal{R}^A{}_B = d\omega^A{}_B - \omega^A{}_D \wedge \omega^D{}_B$]

**D) The (ABAB) Riemann tensor component (RTCs)**

Therefore, adding (7e) and (8c) we obtain::

$$R^{AB}{}_{AB} = +\ [g_A\ g_B]^{-1}\{[g_{A,B}/g_B]_{,B} + \eta_A\ \eta_B[g_{B,A}/g_A]_{,A}\} - \eta_D \eta_B\ g_D^{-2}\ [g_{A,D}/g_A][g_{B,D}/g_B]\ _D \ne _{A,B}$$

….

**Appendix**: **notation**:
 Gauss's curvature for a 2-d metric space (for 2-d one can always find a diagonal metric):
$$ds^2 = E du^2 + G dv^2$$

In terms of the quantities:

$$e = \sqrt{|E|},\ g = \sqrt{|G|} \qquad, \varepsilon \text{ is the sign}$$

$$ds^2 = \varepsilon_1\ e^2 du^2 + \varepsilon_2\ g^2 dv^2$$

Using the notation: $g_u = \partial_u g$ , $(g_u/e)_u = \partial_u(\partial_u g/e)$

the curvature of the surface is given by Gauss's formula:

$$K_{uv} = (-1/eg)[\varepsilon_1(g_u/e)_u + \varepsilon_2(e_v/g)_v]$$

## Appendix: Landau Lifschitz

The Gaussian curvature's essential aspect is: $\left. \dfrac{-g_i}{g_e} \left[ \dfrac{g_{e,i}}{g_i} \right] \right|_{,i}$ . (1a)

The below formula for the Riemann tensor component (RTC) formula for diagonal (orthogonal) metrics, slightly adapted, is from L. D. *Landau*, E. M. *Lifshitz* (section 92)]: where ε contains the signature, and for $i \neq l$:

$R_{lili} = \{\varepsilon_l\, e^{2F}{}_l\, (F_{i,i}F_{l,i} - F^2{}_{l,i} - F_{l,i,i}) + \varepsilon_i\, e^{2F}{}_i(F_{l,l}F_{i,l} - F^2{}_{i,l} - F_{i,l,l})\} - \varepsilon_l\, \varepsilon_i\, \varepsilon_m\, e^{2F}{}_l\, e^{2F}{}_i \Sigma_m \neq {}_{i,l}\, e^{-2F}{}_m\, F_{i,m}\, F_{l,m}$   $i \neq l$   (A)

We do not actually need the { }, however what we wish to show is that in the version with one index raised, this term (the sum of two subterms) is exactly the Gaussian curvature for the two-d represented by (l,i).
Since the two subterms in { } are identical, basically we want to show that:

$\varepsilon_l\, e^{2F}{}_l\, (F_{i,i}F_{l,i} - F^2{}_{l,i} - F_{l,i,i})$ = the Gaussian curvature's essential term (1a).

We demonstrate this in another paper.

## Appendix: "Mathpages"

$$R_{abab} = -\dfrac{g_{aa,bb} + g_{bb,aa}}{2} + \dfrac{1}{4}\left( \dfrac{g_{aa,b}{}^2 + g_{aa,a}g_{bb,a}}{g_{aa}} + \dfrac{g_{bb,a}{}^2 + g_{bb,b}g_{aa,b}}{g_{bb}} - \dfrac{g_{aa,c}g_{bb,c}}{g_{cc}} - \dfrac{g_{aa,d}g_{bb,d}}{g_{dd}} \right)$$

…

## Appendix:
See also: "Ricci Tensor of Diagonal Metric": K. Z. Win

…

### Appendix: Examples of calculation

(employing **the above Formula for the diagonal Riemann Tensor Components** Send feedback to: air1@nyu.edu

Ie using: $R^{AB}{}_{AB} = -[g_A\, g_B]^{-1}\{\eta_A[g_{A,B}/g_B]_{,B} + \eta_B[g_{B,A}/g_A]_{,A}\} - \eta_D\eta_A\eta_B\ g_D{}^{-2}\ [g_{A,D}/g_A][g_{B,D}/g_B]$ $_{D \neq A,B}$

…………..

**Example:** $ds^2 = dt^2 - t^2(dr^2 + \sinh^2 r\, d^2\Omega)$

$g_t = 1, \eta_t = +1;\ g_r = t, \eta_r = -1;\ g_\theta = t \sinh r,\ g_\phi = t \sinh r \sin\theta$.

This certainly seems curved, but all the RTCs vanish, so it is a flat spacetime! Using our formula the calculation of the RTCs is rather simple.

$\mathbf{R^{tr}{}_{tr}} = -[g_t\, g_r]^{-1}\{\eta_t[g_{t,r}/g_r]_{,r} + \eta_r[g_{r,t}/g_t]_{,t}\} - \eta_D\eta_t\eta_r\ g_D{}^{-2}\ [g_{t,D}/g_t][g_{r,D}/g_r]$ $_{D \neq t,r}$

$$= -[t]^{-1}\{ [0/g_r],_r -[1/1],_t\} + \eta_D\ g_D^{-2}\ [0/g_t][g_{r,D}/g_r]\ _{D\neq t,r} = 0$$

$$\mathbf{R^{r\phi}{}_{r\phi}} = -[g_r\ g_\phi]^{-1}\{\eta_r[g_{r,\phi}/g_\phi],_\phi + \eta_\phi\ [g_{\phi,r}/g_r],_r\} - \eta_D\eta_r\eta_\phi\ g_D^{-2}[g_{r,D}/g_r][g_{\phi,D}/g_\phi]\ _{D=t,\theta}$$

In this case, both $g_r$ and $g_\phi$ ARE functions of one of the two intermediate coordinates $(t,\theta)$, namely t, but $g_r$ is not a function of the other intermediate coordinate, $\theta$, and so one of the intermediate terms vanish, ie:

$$= -[t\ g_\phi]^{-1}\{\eta_r[t,_\phi/g_\phi],_\phi + \eta_\phi\ [g_{\phi,r}/g_r],_r\} - \eta_t\ g_t^{-2}[g_{r,t}/g_r][g_{\phi,t}/g_\phi] - \eta_\theta\ g_\theta^{-2}[g_{r,\theta}/g_r][g_{\phi,\theta}/g_\phi]$$

$$= -[t(t\sinh r\sin\theta)]^{-1}\{\eta_r[0/g_\phi],_\phi+\eta_\phi[(t\sinh r\sin\theta),_r/t],_r\}-\eta_t g_t^{-2}[t,_t/t][t\sinh r\sin\theta],_t/[t\sinh r\sin\theta]-\eta_\theta g_\theta^{-2}[t,_\theta/g_r][g_{\phi,\theta}/g_\phi]$$

$$= -[t^2\sinh r\sin\theta)]^{-1}\eta_\phi\ [(t\sinh r\sin\theta),_r/t],_r - \eta_t\ [1]\ [t,_t/t][t,_t/t] - \eta_\theta g_\theta^{-2}\ [0]\ [g_{\phi,\theta}/g_\phi]$$

$$= - t^{-2}\ [\sinh r\sin\theta]^{-1}\eta_\phi\sin\theta\ (\sinh r),_r,_r - \eta_t\ t^{-2}$$

$$= - t^{-2}\{\eta_\phi\ (\sinh r)^{-1}(\sinh r),_r,_r + \eta_t\}$$

$$= - t^{-2}\{\eta_\phi + \eta_t\} = 0.$$

Our intuition is not capable of deciphering a metric directly, looking at it and determining whether or not it represents curvature of flatness. The Riemann tensor provides a fully dependable way to make such a determination: if it vanishes there's no curvature, if it doesn't there IS curvature. It is non-zero if and only if there is non-zero curvature. Whatever the functions $g_A$ and $g_B$ are, however much they make the metric look as though it is that of a curved surface, if **all** components of the Riemann tensor vanishes the space is in fact flat; if even one component is non-zero, it is **NOT** flat.

…

**For angular 2-space in flat 4-d**:

$$\mathbf{R^{r\phi}{}_{r\phi}} = -[1\text{x}\ g_\phi]^{-1}\{\eta_r[1,_\phi/g_\phi],_\phi + \eta_\phi\ [g_{\phi,r}/1],_r\} - \eta_t\ g_t^{-2}[g_{r,t}/g_r][g_{\phi,t}/g_\phi] - \eta_\theta\ g_\theta^{-2}[1,_\theta/1][g_{\phi,\theta}/g_\phi]$$
$$= \eta_\phi\ [1],_r - \eta_t\ g_t^{-2}[1,_t/g_r][g_{\phi,t}/g_\phi] = 0.$$

…

**flat 4-d** in polar coordinates: $ds^2 = -dt^2 + dr^2 + r^2 d\theta^2 + r^2\sin^2\theta d\phi^2$

$$g_r = 1,\ g_t = 1,\ g_\theta = r,\ g_\phi = r\sin\theta$$

$$\mathbf{R^{r\theta}{}_{r\theta}} = -[g_r\ g_\theta]^{-1}\{\eta_r[g_{r,\theta}/g_\theta],_\theta + \eta_\theta\ [g_{\theta,r}/g_r],_r\} - \eta_D\eta_r\eta_\theta\ g_D^{-2}[g_{r,D}/g_r][g_{\theta,D}/g_\theta]\ _{D=t,\phi}$$

$$= -[1\text{x}\ g_\theta]^{-1}\{\eta_r[1,_\theta/g_\theta],_\theta + \eta_\theta\ [g_{\theta,r}/1],_r\} - \eta_t\ g_t^{-2}[g_{r,t}/g_r][g_{\theta,t}/g_\theta] - \eta_\phi\ g_\phi^{-2}[1,_\phi/1][g_{\theta,\phi}/g_\theta]$$
$$= 0 + \eta_\theta\ [1],_r - \eta_t\ g_t^{-2}[1,_t/g_r][g_{\phi,t}/g_\phi] = 0.$$

**.......**

$$\mathbf{R^{r\phi}{}_{r\phi}} = -[g_r\ g_\phi]^{-1}\{\eta_r[g_{r,\phi}/g_\phi],_\phi + \eta_\phi\ [g_{\phi,r}/g_r],_r\} - \eta_D\eta_r\eta_\phi\ g_D^{-2}[g_{r,D}/g_r][g_{\phi,D}/g_\phi]\ _{D=t,\theta}$$

Since $g_r$ is not a function of $\phi$, obviously the first part of the { } will vanish, and since $g_r$ is not a function of either of the two intermediate coordinates $(t,\theta)$, the last two terms vanish, ie:

$$= -[1 \times g_\phi]^{-1}\{\eta_r[1,_\phi/g_\phi],_\phi + \eta_\phi [g_{\phi,r}/1],_r\} - \eta_t\, g_t^{-2}[g_{r,t}/g_r][g_{\phi,t}/g_\phi] - \eta_\theta\, g_\theta^{-2}[g_{r,\theta}/g_r][g_{\phi,\theta}/g_\phi]$$

$$= -g_\phi^{-1}\eta_\phi\, g_{\phi,r,r} - \eta_t\, g_t^{-2}[0/g_r][0/g_\phi] - \eta_\theta\, g_\theta^{-2}[0/g_r][g_{\phi,\theta}/g_\phi]$$

$$= -[r\sin\theta]^{-1}\eta_\phi\, [r\sin\theta],_{r,r} \;,$$

however: $[r\sin\theta],_{r,r} = [\sin\theta]\, r,_{r,r} = 0$, so $R^{r\phi}{}_{r\phi} = 0$.

After some experience with this type of calculation, just by inspection one can write:

$$\mathbf{R^{r\phi}{}_{r\phi}} = -[r\sin\theta]^{-1}\eta_\phi\, [r\sin\theta],_{r,r} = 0$$

…………..

$$\mathbf{R^{\theta\phi}{}_{\theta\phi}} = -[g_\theta\, g_\phi]^{-1}\{\eta_r[g_{\theta,\phi}/g_\phi],_\phi + \eta_\theta[g_{\phi,\theta}/g_\theta],_\theta\} - \eta_t\eta_\theta\eta_\phi\, g_t^{-2}[g_{\theta,t}/g_\theta][g_{\phi,t}/g_\phi] - \eta_r\eta_\theta\eta_\phi\, g_r^{-2}[g_{\theta,r}/g_\theta][g_{\phi,r}/g_\phi]$$

$$= -[g_\theta\, g_\phi]^{-1}\{\eta_r[0/g_\phi],_\phi + \eta_\theta[r\cos\theta/r],_\theta\} - \eta_t\eta_\theta\eta_\phi\, g_t^{-2}[0/g_\theta][g_{\phi,t}/g_\phi] - \eta_r\eta_\theta\eta_\phi\, 1^{-2}[r,_r/g_\theta][(r\sin\theta),_r/g_\phi]$$

$$= -[g_\theta\, g_\phi]^{-1}\eta_\theta\, [-\sin\theta]\} - \eta_r\eta_\theta\eta_\phi\, [1/g_\theta][\sin\theta/g_\phi] = 0$$

**Analysis**: the cancellation is from: $-[g_\theta\, g_\phi]^{-1}\eta_\theta\, [g_{\phi,\theta}/g_\theta],_\theta = \eta_r\eta_\theta\eta_\phi\, g_r^{-2}[g_{\theta,r}/g_\theta][g_{\phi,r}/g_\phi]$

$- [g_{\phi,\theta}/g_\theta],_\theta = \eta_r\eta_\phi\, g_r^{-2}\, g_{\theta,r}\, g_{\phi,r}$

$- [g_{\phi,\theta}/g_\theta],_\theta/[\eta_\phi\, g_{\theta,r}\, g_{\phi,r}] = \eta_r\, g_r^{-2}$ : for $\eta_r\, g_r = 1$, this is: $[g_{\phi,\theta}/g_\theta],_\theta = -[\eta_\phi\, g_{\theta,r}\, g_{\phi,r}]$

For $g_\theta = r$ : $[g_{\phi,\theta}/r],_\theta = -[\eta_\phi\, r,_r\, g_{\phi,r}]$ , $[1/r]g_{\phi,\theta,\theta} = -\eta_\phi\, g_{\phi,r}$

For $\eta_\phi = 1$, $g_\phi = rf(\theta)$ : $f(\theta),_{\theta,\theta} = -[rf(\theta)],_r = -f(\theta)$ and so $f = \sin$ (maybe cos also?)